\documentclass[a4paper]{article}

\let\ucal=\mathcal

\usepackage{amsmath}
\usepackage{amsthm}
\usepackage{amssymb}
\usepackage{eucal}
\newtheorem{pro}{Proposition} 
\newtheorem{lem}[pro]{Lemma}
\newtheorem{thm}[pro]{Theorem}
\newtheorem{corl}[pro]{Corollary}
\theoremstyle{remark}
\newtheorem{rem}[pro]{Remark}
\theoremstyle{definition}

\numberwithin{equation}{section}

\def\toRomans{
  \renewcommand{\theenumi}{(\roman{enumi})}}

\def\toAlpha{
  \renewcommand{\theenumi}{(\alph{enumi})}}

\toRomans

\def\and{and }

\title{The Haar system in the preduals of hyperfinite factors.}
\author{D.~Potapov\thanks{Corresponding author, e-mail:
    pota0002@infoeng.flinders.edu.au \endgraf \quad The research is
    partially supported by ARC}\ \ \and F.~Sukochev} 

\def\HFF{{\mathcal{R}}}
\def\aA{{\mathcal{A}}}

\def\Bd{B}
\def\cC{{\ucal{C}}}
\def\aM{{\mathcal{M}}}
\def\aN{{\mathcal{N}}}
\def\sL{{\mathcal{L}}}
\def\Q{{\mathbb{Q}}}
\def\N{{\mathbb{N}}}

\def\sH{{\mathcal{H}}}
\def\Rl{{\mathbb{R}}}
\def\Cx{{\mathbb{C}}}

\def\id{{\mathbf{1}}}

\def\dtimes{{\times\kern-1pt\times}}

\def\e{{\mathbf{e}}}
\def\h{{\mathbf{h}}}

 \def\x{{\mathbf{x}}} \def\y{{\mathbf{y}}}
\def\z{{\mathbf{z}}}

\def\D{{\mathbf{D}}}
\def\H{{\mathbf{H}}}
\def\E{{\mathbf{E}}}

\begin{document}

\bibliographystyle{short}

\baselineskip=1.2\baselineskip
\parskip=5pt plus 2pt minus 1pt
\tolerance=400


\maketitle

\begin{abstract}
We shall present examples of Schauder bases in the preduals to the
hyperfinite factors of types~$\hbox{II}_1$, $\hbox{II}_\infty$,
$\hbox{III}_\lambda$, $0 < \lambda \leq 1$.  In the semifinite
(respectively, purely infinite) setting, these systems form Schauder bases
in any associated separable symmetric space of measurable operators
(respectively, in any non-commutative $L^p$-space).
\end{abstract}

\subsubsection{Introduction}

A sequence~$\x = \{x_n\}_{n \geq 1}$ in a Banach space~$X$ is called a
(Schauder) basis of~$X$ if, for every~$x\in X$ there exists a unique
sequence of scalars~$\{\alpha_n\}_{n \geq 1}$ so that $$ x = \sum_{n
  \geq 1} \alpha_n x_n. $$ A sequence~$\x$ such that~$x_n \ne 0$ for
all~$n$ and~$[x_n]_{n \geq 1} = X$ forms a basis of~$X$ if and only if
there is a constant~$c$ so that for every choice of
scalars~$\{\alpha_n\}_{j=1}^k$ and integers~$m < k$ we have $$
\|\sum_{1 \leq j \leq m} \alpha_j x_j \|_X \le c\, \|\sum_{1 \leq j
  \leq k} \alpha_j x_j \|_X\ \ \hbox{(cf.~\cite{LT-I}).} $$ The smallest
such constant~$c$ is called the basis constant of~$\x$.  In this note
we shall be concerned with the construction of Schauder bases in
spaces of operators associated with the hyperfinite factors of
type~$\hbox{II}$ and~$\hbox{III}_\lambda$, $0 < \lambda \leq 1$.  In
the setting of symmetric spaces of measurable operators affiliated
with the hyperfinite factors of type~$\hbox{II}$ and with some
hyperfinite von Neumann algebras of type~$\hbox{I}_\infty$ the problem
was recently considered in~\cite{SF1994, SF1995, DFPS2001} where
``non-commutative Walsh system'', ``non-commutative trigonometric
system'' and ``non-commutative Vilenkin systems'' were constructed.
However, as with their classical counterparts, these systems fail to
form a Schauder basis in the preduals.  In order to construct a
Schauder basis in the preduals to the hyperfinite factors, we use an
analogy with another classical function system (which forms a Schauder
basis in every separable symmetric function space on~$(0, 1)$
\cite{KPS-IOLO, LT-I}), namely with the Haar system.

{\bf Acknowledgements.} We would like to thank professor J.~Arazy and
professor A.~Pe{\l}czy{\'n}ski for useful discussions.

\subsubsection{Preliminaries.}

\def\CE{{\mathcal{E}}}

Let~$\aM$ be a von Neumann algebra with a fixed faithful normal
state~$\rho$.  Let~$\aM_*$ stand for the predual of~$\aM$, we consider
this predual as a subspace of~$\aM^*$ consisting of all normal linear
functionals equipped with~$\|\cdot\|_* := \|\cdot\|_{\aM^*}$,
cf.~\cite[Theorem~1.10]{SZ-Lec-on-vNA}.

We consider several different norms on~$\aM$.  $\|\cdot\|$ is the
operator norm.  The norm~$\|\cdot\|_*$ of the predual~$\aM_*$ induces
the norms~$\|\cdot\|_\sharp$ and~$\|\cdot\|_\flat$ on~$\aM$ by means of
the left and right embeddings~$x\in\aM \rightarrow x\rho:=\rho(\cdot
x)\in\aM_*$ and~$x\in\aM \rightarrow \rho x := \rho(x \cdot)\in\aM_*$,
respectively.  These embeddings are injective and the ranges of~$\aM$
under these embeddings are dense in~$\aM_*$, \cite{Ko1984, Ray2003}.
Thus, if~$\aM_\sharp$ and~$\aM_\flat$ are completions of~$\aM$ with respect
to the norms~$\|\cdot\|_\sharp$ and~$\|\cdot\|_\flat$, then these spaces
are isometric to~$\aM_*$.  Obviously,~$\|\cdot\| \ge \|\cdot\|_\sharp$
and~$\|\cdot\| \geq \|\cdot\|_\flat$, and the embeddings~$\aM
\subseteq \aM_\sharp$ and~$\aM \subseteq \aM_\flat$ are continuous.
Let us also note that the space~$\aM_{\sharp}$ (resp.~$\aM_\flat$) is
a left (resp.\ right) module with respect to~$\aM$, i.e.\ $xa \in
\aM_\sharp$ (resp.\ $ax \in \aM_{\flat}$) provided~$a \in \aM$ and~$x
\in \aM_{\sharp}$ (resp.~$x \in \aM_{\flat}$), moreover we have $$
\|xa\|_{\aM_\sharp} \leq \|x\|_{\aM_\sharp}\, \|a\|\ \ (\text{resp.}\
\|ax\|_{\aM_\flat} \leq \|a\|\, \|x\|_{\aM_\flat}). $$

We introduce left (resp.\ right) $L^p$-space associated with the
algebra~$\aM$ as
\begin{equation}
  \label{LpInterp}
  L^p_{\sharp (\text{resp.}\  \flat)}(\aM) := [\aM, \aM_{\sharp
  (\text{resp.}\    \flat)}]_{\frac 1p},\ \ 1\le p\le \infty. 
\end{equation}
Here,~$[\cdot, \cdot]_\theta$ is the method of complex
interpolation,~\cite{BergLofstrom}.  The space~$L^p_{\sharp
  (\text{resp.}\ \flat)}(\aM)$ is isomorphic to Haagerup's
$L^p$-spaces $L^p(\aM)$, \cite{Terp1981}.  Clearly, $L^1_{\sharp
  (\text{resp.}\ , \flat)}(\aM) = \aM_{\sharp (\text{resp.}\ \flat)}$
and~$L^\infty_\sharp = L^\infty_\flat = \aM$.  Moreover, the Hilbert
space~$L^2_{\sharp (\text{resp.}\ \flat)}(\aM)$ coincides with the
completion of~$\aM$ with respect to the inner product~$\langle x, y
\rangle_\sharp := \rho(y^*x)$ (resp.\ $\langle x, y \rangle_\flat :=
\rho(xy^*)$), $x, y \in \aM$.  We refer the reader to~\cite{Ko1984,
  Ray2003}, for further details on this construction and also
to~\cite{Terp1981} for construction of Haagerup's $L^p$-spaces.

In this note, we shall prove the results for the left
norm~$\|\cdot\|_\sharp$.  The argument for the right
norm~$\|\cdot\|_\flat$ is generally the same.  We shall make appropriate
remarks when it is necessary.

We denote by~$\sigma^\rho$ the modular automorphism group for the
state~$\rho$, i.e. the unique automorphism group of~$\aM$ such that
(i)~$\rho(x) = \rho(\sigma^\rho_t(x))$, $t\in\Rl$, $x\in\aM$ and~(ii)
for every~$x, y\in\aM$, there is a complex function~$f_{x,y}(z)$
bounded in the strip~$\bar S$ and holomorphic in~$S$, where~$S = \{z
\in \Cx:\ \ 0< \hbox{Im}\ z < 1 \}$ such that~$\rho(\sigma^\rho_t(x)\,
y) = f_{x, y}(t)$ and~$\rho(y\, \sigma^\rho_t(x)) = f_{x, y}(t + i)$,
cf.~\cite[Section~9.2]{KadRingII}.

Let~$\aN \subseteq \aM$ be a von Neumann subalgebra such that the
modular group~$\sigma^\rho$ (with respect to~$\aM$) leaves~$\aN$
globally invariant.  In this case the restriction of~$\sigma^\rho$
onto~$\aN$ gives the modular group of~$\rho|_\aN$ in the algebra~$\aN$
and we can speak about the modular action~$\sigma^{\rho}$ without
referring to the particular algebra~$\aN$ or~$\aM$.  In this (and only in
this) setting, according to the main result of~\cite{Ta1972}, there is a
norm one projection~$\CE : \aM \rightarrow \aN$, which we call
conditional expectation, such that \toAlpha
\begin{enumerate}
\item $\rho(x) = \rho(\CE x)$, $x\in\aM$; \label{EIprop}
\item $\CE(a x b) = a \CE(x) b$, $a, b\in \aN$, $x\in\aM$;
\item $0 \le \CE(x)^* \CE(x) \le \CE(x^* x)$.\label{EIIIprop}
\end{enumerate}
\toRomans

We now show that the existence of the conditional expectation
implies that (i)~$\aN_\sharp$ continuously embeds into~$\aM_\sharp$ and
(ii)~the space~$\aN_{\sharp}$ is a $1$-complemented subspace
in~$\aM_{\sharp}$.  (i)~follows from the inequality
\begin{equation}
  \label{ContEmbed}
  \begin{aligned}
    \|x\|_{\aM_\sharp} := &\, \|x \rho\|_{\aM_*} = \sup_{y \in \aM}
    \|y\|^{-1}\, | \rho(y x)| \\ \text{(since~$x\in\aN$)}\ \ = &
    \sup_{y\in\aM} \|y\|^{-1}\, | \rho(\CE(y) x)| \\
    \text{(since~$\|\CE(y)\|\le \|y\|$)}\ \ \le &\, \sup_{y\in\aN}
    \|y\|^{-1} \, |\rho(y x)| \\ = &\, \|x\rho\|_{\aN_*} =
    \|x\|_{\aN_\sharp},\ \ x\in \aN.
  \end{aligned}
\end{equation}
For~(ii), it is sufficient to show that
\begin{equation}
\label{bddPrj}
\|\CE(x)\|_{\aN_{\sharp}} \le
\|x\|_{\aM_{\sharp}},\ \ x\in\aM. 
\end{equation}
Let us recall, that the predual~$\aM_*$ is a subspace of~$\aM^*$
consisting of all normal linear functionals.  Let us consider the
mapping~$\CE' : \aM_* \rightarrow \aN_*$ given by~$\CE'(\phi) =
\phi|_{\aN}$.  It follows from the
properties~\ref{EIprop}--\ref{EIIIprop} of~$\CE$ above that
\begin{equation}
  \label{PrjIdentity}
  \CE'(x\rho) = \CE(x) \rho\ \ \text{and}\ \ \CE'(\rho x) = \rho
  \CE(x),\ \ x\in\aM.
\end{equation}
Now~\eqref{bddPrj} follows from~\eqref{PrjIdentity}, since $\CE'$ is a
norm one linear operator.  It follows from~\eqref{ContEmbed}
and~\eqref{bddPrj} that the embedding~$\aN_\sharp \subseteq
\aM_\sharp$ is isometric and therefore, the space~$\aN_\sharp$ is a
$1$-complemented subspace of~$\aM_\sharp$.

It this paper, we only consider von Neumann subalgebras~$\aN \subseteq
\aM$, which are globally invariant under~$\sigma^\rho$.  Therefore, we
shall refer to the norms in the spaces~$\aN_{\sharp}$
and~$\aM_{\sharp}$ simply as~$\|\cdot\|_{\sharp}$ without specifying
the particular algebra.

The assertions above may be similarly carried to the right predual space
and to~$L^p$ spaces (left and right) by interpolation.  That is, the
space~$L^p_{\sharp (\text{resp.}\ \flat)}(\aN)$ is $1$-complemented
in~$L^p_{\sharp (\text{resp.}\ \flat)} (\aN)$, $1\le p \le \infty$ and
$L^p_{\sharp (\text{resp.}\ \flat)}(\aN)$ embeds isometrically
into~$L^p_{\sharp (\text{resp.}\ \flat)}(\aM)$, $1 \leq p \leq \infty$.

\subsubsection{Matrix spaces.}

\def\Alpha{\hbox{A}}

Let $\nu\in\N$ and~$0<\alpha \le \frac 12$ be fixed
throughout the text.  Let~$\aN_\nu$ be the class of all complex $2^\nu
\times 2^\nu$-matrices with the unit matrix~$\id_\nu$.  $Tr$ is the standard
trace on matrices.  The state~$\rho_\nu$ on~$\aN_\nu$ is given by
\begin{equation}
  \label{statedef}
  \rho_\nu(x) = Tr (x \Alpha_\nu),\ \ x\in\aN_\nu,\ \ \Alpha_\nu =
  \bigotimes_{k=1}^\nu \left[\begin{matrix} \alpha & 0 \\ 0 & 1 -
  \alpha \\   \end{matrix} \right].
\end{equation}
The definition of the state~$\rho_\nu$ immediately implies that
\begin{equation}
  \label{tensorProp}
  \rho_{\nu+\mu} (x \otimes y) = \rho_\nu (x)\, \rho_\mu(y),\ \
  x\in\aN_\nu, y\in\aN_\mu.
\end{equation}
We consider the ultra-weak continuous $*$-isomorphic embedding~$i_\nu:
\aN_\nu \rightarrow \aN_{\nu+1}$ given by
\begin{equation}
  \label{TowerEmbed}
  i_\nu(x) = x \otimes \id_1,\ \ x\in\aN_\nu.
\end{equation}
Due to~\eqref{tensorProp}, we have that~$\rho_{\nu+1}(i_\nu(x)) =
\rho_\nu(x)$, $x\in\aN_\nu$, i.e. the restriction of the
state~$\rho_{\nu+1}$ onto the subalgebra~$i_\nu(\aN_\nu)$ is equal to the
state~$\rho_\nu$.

The modular automorphism group of the state~$\rho_\nu$ is given
by~$\sigma^{\rho_{\nu}}_t(x) = \Alpha_\nu^{it} x \Alpha^{-it}_\nu$,
$t\in\Rl$, $x\in\aN_\nu$.  Indeed,
\begin{align*}
\rho_\nu(\sigma_t^{\rho_\nu}(x)\, y) = &\, Tr (\Alpha_\nu^{1+it}
x \Alpha_\nu^{-it} y) = f_{x, y}(t), \\
\rho_\nu(y\, \sigma^{\rho_\nu}_t(x)) = &\, Tr (\Alpha_\nu^{it} x
\Alpha^{1-it}_k y) = f_{x, y}(t + i), 
\end{align*}
where the holomorphic function~$f_{x, y}$ is given by~$f_{x,y}(z) = Tr
(\Alpha_\nu^{1+iz} x \Alpha_\nu^{-iz} y)$, $x, y\in \aN_\nu$.
Therefore, it is readily seen that the group~$\sigma^{\rho_{\nu+1}}$ leaves the
subalgebra~$i_\nu(\aN_\nu)$ globally invariant.  According to the
preceding section, the space~$i_\nu(\aN_{\nu, \sharp})$ is
$1$-complemented in~$\aN_{\nu+1, \sharp}$ and the mapping~$i_\nu$
embeds~$\aN_{\nu, \sharp}$ isometrically into $\aN_{\nu+1, \sharp}$.  We
denote~$\CE_\nu$ the norm one projection~$\aN_{\nu+1, \sharp}
\rightarrow i_\nu (\aN_{\nu, \sharp})$.  From now on, we shall refer to
the norms in the spaces~$\aN_{\nu, \sharp}$ simply
as~$\|\cdot\|_{\sharp}$ omitting the index~$\nu$.  

Similarly, we introduce the $L^p$-space~$\sL^p_{\nu, \sharp} :=
L^p_\sharp (\aN_\nu)$, $\nu \geq 1$ and we refer to the norm in this
space as~$\|\cdot\|_{\sharp , p}$.  

We also introduce the $p$-th Schatten-von Neumann
norm~$\|\cdot\|_{\cC_p}$, $1\leq p < \infty$ on~$\aN_\nu$ as $$
\|x\|_{\cC_p} := (Tr((x^*x)^{\frac p 2}))^{\frac 1p},\ \ x\in\aN_\nu. $$
$\|\cdot\|_{\cC_\infty}$ stands for the operator norm.  We
denote~$\cC^{(\nu)}_p$ the matrix space~$\aN_\nu$ equipped with the
norm~$\|\cdot\|_{\cC_p}$, $1 \leq p \leq \infty$.  We now may express
the norms~$\|\cdot\|_\sharp$ and~$\|\cdot\|_\flat$ as
\begin{equation}
  \label{leftRightPreNorms}
  \|x\|_{\sharp} = \|x \Alpha_\nu\|_{\cC_1},\ \ \|x\|_{\flat} =
  \|\Alpha_\nu x \|_{\cC_1},\ \ x\in\aN_\nu.
\end{equation}
The last identities may be carried to the $L^p$-spaces associated
with~$\aN_\nu$ as follows.

\begin{rem}
  \label{LpIdent}
  We fix~$\nu\in\N$ and consider the function~$f^\sharp: \Cx \times
  (\aN_\nu + \aN_{\nu \sharp}) \rightarrow \cC^{(\nu)}_\infty +
  \cC^{(\nu)}_1$ given by $$ f^\sharp (z, x) = x \Alpha_\nu^z,\ \
  z\in\Cx,\ \ x\in\aN_\nu + \aN_{\nu, \sharp}. $$ For every
  fixed~$z\in\Cx$, $f^\sharp _z (\cdot) := f^\sharp (z, \cdot)$ is a
  linear operator~$\aN_\nu + \aN_{\nu, \sharp} \rightarrow
  \cC^{(\nu)}_\infty + \cC^{(\nu)}_1$.  Thus, we may
  consider~$f^\sharp_{(\cdot)}$ as a function on the complex plane
  with values in~$B(\aN_\nu + \aN_{\nu, \sharp}, \cC^{(\nu)}_\infty +
  \cC^{(\nu)}_1)$.  Here~$B(X, Y)$ is the Banach space of all bounded
  linear operators~$X \rightarrow Y$.  The
  function~$f^\sharp_{(\cdot)}$ is holomorphic on~$0 < \hbox{Re}\ z <
  1$.  It follows from~\eqref{leftRightPreNorms} that the
  mapping~$f^{\sharp}_{1 + it}$ is an isometry between~$\cC^{(\nu)}_1$
  and~$\aN_{\nu, \sharp}$, for every~$t\in\Rl$.  On the other hand,
  the mapping~$f^{\sharp}_{it}$ is clearly an isometry
  between~$\cC^{(\nu)}_\infty$ and~$\aN_\nu$, for every~$t\in\Rl$.
  Thus, interpolating, we obtain that the mapping~$f^\sharp_{\frac
    1p}$ is an isometry between~$\cC^{(\nu)}_p$ and~$\sL^p_{\nu, \sharp}$, i.e.
  \begin{equation}
    \label{LpLeftRightIdent}
    \|x\|_{p, \sharp} = \|x \Alpha^{\frac 1p}_\nu\|_{\cC_p},\ \
    x\in\sL^p_{\nu, 
      \sharp}. 
  \end{equation}
  Similar argument for the right spaces gives
  \begin{equation}
    \label{LpRightIdent}
  \|x\|_{p, \flat} = \|\Alpha^{\frac 1p}_\nu x\|_{\cC_p},\ \ x\in\sL^p_{\nu,
  \flat}. 
  \end{equation}
\end{rem}

Let~$\e_\nu = \{e^{(\nu)}_j\}_{0\le j< 4^\nu}$ be the matrix units
in~$\aN_\nu$ given in the shell enumeration,
cf.~\cite{KwPe1970,Ar1978}, see also the proof of
Theorem~\eqref{TesorBasis}.  The system of matrix units~$\e_\nu$ forms
a basis of~$\cC_p$ with the basis constant~$2$,
cf.~\cite[Theorem~2.1]{KwPe1970}.  The next lemma shows that this
system remains a basis with the same basis constant~$2$ with respect
to the norms~$\|\cdot\|_{p, \sharp}$ and~$\|\cdot\|_{p, \flat}$, $1\le
p \le \infty$.

\begin{lem}
  \label{MatrixUnits} For every~$0\le m < 4^\nu$, every~$1\le p \le
  \infty$ and any complex numbers~$\alpha_j \in \Cx$, $0\le j <
  4^\nu$, we have $$ \|\sum_{0\le j \le m} \alpha_j\, e^{(\nu)}_j \|_{p,
    \sharp (\text{resp.}\ \flat)} \le 2\, \|\sum_{0\le j < 4^\nu}
  \alpha_j e^{(\nu)}_j \|_{p, \sharp (\text{resp.}\ \flat)}. $$
\end{lem}

\begin{proof} 
  Let~$P_{\nu, m}$ be the basis projection corresponding to the
  number~$0\le m <4^\nu$.  The projection~$P_{\nu, m}$ is a Schur
  multiplier, i.e.\ 
  $$
  P_{\nu, m} (x) = p_{\nu, m} \circ x, $$
  where~$\circ$ is the Schur
  (entrywise) product of matrices and $$
  p_{\nu, m} = \sum_{0 \leq j <
    m} e^{(\nu)}_j. $$
  Let us note that the Schur product is commutative
  and multiplication by a diagonal matrix is a special case of Schur
  multiplier.  Thus, the claim of the lemma follows from the
  result~\cite[Theorem~2.1]{KwPe1970}, the
  identities~\eqref{leftRightPreNorms}, \eqref{LpLeftRightIdent},
  \eqref{LpRightIdent} and the fact that the operator~$P_{\nu , m}$
  commutes with left and right multiplication by a diagonal matrix, i.e.
  $$
  P_{\nu, m}(x \Alpha^{1/p}_\nu) = P_{\nu, m}(x)\, \Alpha^{1/p}_\nu\ 
  \ \text{and}\ \ P_{\nu, m}(\Alpha^{1/p}_\nu x) = \Alpha^{1/p}_\nu \,
  P_{\nu, m}(x).
  $$
\end{proof}

At the end of the section we establish the explicit formula of the
projection~$\CE_\nu$ on elementary tensors, i.e.
\begin{equation}
  \label{CEprop}
  \CE_\nu(x \otimes y) = \rho_1(y)\, i_\nu(x),\ \ x\in\aN_\nu, y\in\aN_1. 
\end{equation}
To this end, consider the Hilbert spaces~$\sH^\sharp_\nu :=
\sL^2_{\nu, \sharp}$, which is the matrix space~$\aN_{\nu}$, equipped
with the inner product~$\langle x, y \rangle_\nu = \rho_\nu(y^*x)$,
$x,y\in\aN_\nu$ and observe that the projection~$\CE_\nu$ as an
orthogonal projection in the Hilbert space~$\sH^\sharp_{\nu+1}$ onto
the subspace~$i_\nu(\sH^\sharp_\nu)$.  If~$\{f_j\}_{0\le j < 4^\nu}$
is an orthonormal basis in~$\sH^\sharp_\nu$, then,
using~\eqref{tensorProp}, we obtain~\eqref{CEprop} as follows
\begin{align*}
  \CE_\nu(x \otimes y) = &\, \sum_{0\le j < 4^\nu} \langle i_\nu(f_j), x
  \otimes y \rangle_{\nu+1}\, i_\nu(f_j) \\ = &\, \sum_{0\le j < 4^\nu}
  \rho_1(y)\, \langle f_j, x \rangle_{\nu}\, i_\nu(f_j)\\ = &\,
  \rho_1(y)\, i_\nu(x),\ \ x \in\aN_\nu, y\in\aN_1.
\end{align*}

\subsubsection{The Haar system.}

We shall construct the Haar system on~$\aN_\nu$ with respect
to~$\rho_\nu$ inductively.  Let us first note that we construct two
Haar systems: the left and the right Haar system, which coincide
when~$\rho_\nu$ is a tracial state.  We shall show the construction of
the left Haar system.  At the outset, we fix an orthonormal basis
of~$\sL^2_{1, \sharp}$, that is elements~$r_j\in\aN_1$, $0\le j \le 3$
such that
\begin{equation}
  \label{basis}
  \rho_1(r^*_j r_k) = \delta_{jk},\ \ 0\le j,k\le 3.
\end{equation}
We define the Haar system inductively.  The Haar system~$\h_1$
in~$\aN_1$ is the system $\h_1 = \{r_0, r_1, r_2, r_3\}$.  If~$\h_\nu =
\{h^{(\nu)}_j\}_{0\le j < 4^\nu}$ is the Haar system in~$\aN_\nu$, then
the system~$\h_{\nu + 1} = \{ h^{(\nu + 1)}_j\}_{0\le j < 4^{\nu+1}}$
given by
\begin{equation}
  \label{Inductive}
  h^{(\nu+1)}_j = 
  \begin{cases} 
    i_\nu(h^{(\nu)}_k) \cdot (\id_\nu
    \otimes r_0), & \text{if~$q = 0$;} \\
    i_\nu(e^{(\nu)}_k) \cdot (\id_\nu \otimes r_q), &
    \text{if~$q\neq 0$;}
  \end{cases},\ \ 0\le j < 4^{\nu+1},
\end{equation}
$$ j = 4^\nu q + k,\ \ 0\le q \le 3,\ \ 0\le k < 4^\nu $$
is the Haar system in~$\aN_{\nu + 1}$.

We shall now present an inductive estimate of the basis constant of
the system~$\h_\nu$ in the space~$\aN_{\nu, \sharp}$.

\begin{thm}
  \label{MTh} If~$c_{\nu, \sharp}$ is the basis constant for~$\h_\nu$,
  then, we have $$ c_{1, \sharp} \le \sum_{0\le j \le 3} \|r_j\|\,
  \|r_j\|_{\sharp} $$ and $$ c_{\nu+1, \sharp} \le \max\{ c_{\nu,
    \sharp}\, \|r_0\|^2, \|r_0\|^2 + 2\, \sum_{q=1}^3 \|r_q\|^2\}. $$
\end{thm}

\begin{proof}
  For the first inequality, it is sufficient to note that, if~$x =
  \sum_{j=0}^3 \alpha_j r_j$, then~$\alpha_j = \rho_1(r^*_j x)$, $0\le
  j\le 3$, see~\eqref{basis}, therefore, $|\alpha_j| \le \|r_j\|\,
  \|x\|_{\sharp}$ and
  $$
  \| \sum_{0\le j \leq m} \alpha_j r_j \|_\sharp \le \sum_{0\le j \leq 3}
  |\alpha_j|\, \|r_j\|_\sharp \le \|x\|_\sharp \sum_{0\le j \le 3} \|r_j\|\,
  \|r_j\|_\sharp,\ \ 0 \leq m < 4. $$

  Let~$r_{j, \nu}:=\id_\nu \otimes r_j$, $0\le j \le 3$.  We shall
  estimate the constant~$c_{\nu+1, \sharp}$ in the inequality
  \begin{equation}
    \label{NextBasisIneq}
    \|\sum_{0\le j \leq m} \alpha_j\, h^{({\nu+1})}_j \|_{\sharp}
    \le c_{{\nu+1}, \sharp}\, \|\sum_{0\le j< 4^{\nu+1}} \alpha_j\,
    h^{({\nu+1})}_j \|_{\sharp},
  \end{equation}
  where~$0\le m < 4^{\nu+1}$.  We first establish the estimate 
  \begin{equation}
    \label{intI}
    \|\sum_{0\le j < 4^\nu} \alpha_j\, h^{(\nu)}_j \|_{\sharp} \le 
    \|r_0\|\, \|\sum_{0\le j < 4^{\nu+1}} \alpha_j\,
    h^{(\nu+1)}_j \|_{\sharp}.
  \end{equation}
  This inequality follows from the observations that
  (cf.~\eqref{CEprop}, \eqref{basis} and~\eqref{Inductive}) $$
  i^{-1}_\nu(\CE_\nu (r_{0, \nu}^*\, h^{(\nu +1)}_j)) =
  \begin{cases} h^{(\nu)}_j,& \text{if~$0\le j < 4^\nu$;} \\ 0, &
    \text{if~$j\ge 4^\nu$}; \end{cases} $$
  and that the left
  multiplication is a bounded operation in~$\|\cdot\|_{\sharp}$.  It
  follows from~\eqref{intI} that, for~$0 \leq m < 4^\nu$, we have
  \begin{align*}
    \|\sum_{0\le j \leq m} \alpha_j h^{(\nu+1)}_j \|_{\sharp} \le &\,
    \|r_0\|\, \| \sum_{0 \leq j \leq m} \alpha_j \, i_\nu(h^{(\nu)}_j)
    \|_\sharp \\ = &\, \|r_0\|\, \|\sum_{0 \leq j \leq m} \alpha_j
    h^{(\nu)}_j \|_\sharp \leq 
    c_{\nu, \sharp}\, \|r_0\|\, \|\sum_{0\le j < 4^\nu} \alpha_j h^{(\nu)}_j
    \|_{\sharp} \\ \le &\, c_{\nu, \sharp} \, \|r_0\|^2\, 
    \|\sum_{0\le j < 4^{\nu+1}} \alpha_j h^{(\nu+1)}_j \|_{\sharp}. \\
  \end{align*}
  
  Therefore, if~$0\le m < 4^\nu$, then the constant
  in~\eqref{NextBasisIneq}, admits the estimate
  \begin{equation}
    \label{resultI}
    c_{\nu+1, \sharp} \le c_{\nu, \sharp}\, \|r_0\|^2.
  \end{equation}
  Let us next establish the estimate
  \begin{equation}
    \label{intII}
    \|\sum_{q 4^\nu \le j \leq q 4^\nu + m} \alpha_j\, h^{(\nu +1)}_j
    \|_{\sharp} \le 2\, \|r_q\|^2\, \|\sum_{0\le j < 4^{\nu+1}}
    \alpha_j\, h^{(\nu+1)}_j \|_{\sharp},
  \end{equation}
  $$
  1\le q \le 3,\ \ 0\le m < 4^\nu. $$
  To this end, we observe that,
  according to~\eqref{CEprop} and~\eqref{basis}, for~$1\leq q \leq 3$ $$
  i^{-1}_\nu(\CE_\nu(r^*_{q,\nu}\, h^{(\nu+1)}_j)) =
  \begin{cases} e^{(\nu)}_{j-q 4^\nu}, & \text{if~$q4^\nu \le j <
      (q+1)4^\nu$}; \\ 0, & \text{otherwise.} \\ \end{cases} $$
  Thus,
  if~$P_{\nu, m}$ is the projection from Lemma~\ref{MatrixUnits}, then
  the left side can be obtained from the right side in~\eqref{intII} via
  the mapping~$x \rightarrow r_{q,\nu}\, i_\nu(P_{\nu,
    m}i^{-1}_\nu(\CE_\nu(r^*_{q,\nu}\, x)))$ and, therefore,
  \eqref{intII} follows.
  
  Finally, if~$4^\nu \le m < 4^{\nu+1}$, then, combining~\eqref{intI}
  and~\eqref{intII}, we obtain for the constant
  in~\eqref{NextBasisIneq} the estimate
  \begin{equation}
    \label{resultII}
    c_{\nu+1, \sharp} \le \|r_0\|^2 + 2\, \sum_{q=1}^3 \|r_q\|^2.
  \end{equation}
  The theorem is proved.
\end{proof}

\begin{rem}
  \label{rightremark}
  For the right Haar system the construction is the same,
  except we start with the system~$\{r_0, r_1, r_2, r_3\}$ such that 
  \begin{equation}
    \label{leftBasis}
    \rho_1(r_j r^*_k) = \delta_{jk},\ \ 0\le j,k\le 3. 
  \end{equation}
  Clearly, in the proof of Theorem~\ref{MTh}, all references to the left
  multiplication should be replaced with those to the right multiplication.
  Thus, we obtain that the right Haar system basis constant~$c_{\nu,
    \flat}$ admits the similar estimates: $$
  c_{1, \flat} \le \sum_{0\le
    j \le 3} \|r_j\|\, \|r_j\|_{\flat} $$
  and $$
  c_{\nu+1, \flat} \le
  \max\{ c_{\nu, \flat}\, \|r_0\|^2, \|r_0\|^2 + 2\, \sum_{q=1}^3
  \|r_q\|^2\}. $$
\end{rem}

\begin{rem}
  \label{LpMth}
  Inspection of the proof of Theorem~\ref{MTh} shows that the main
  ingredient of the proof is (i)~Lemma~\ref{MatrixUnits} and (ii)~the
  fact that the left multiplication by a bounded operator is
  continuous in the norm~$\|\cdot\|_{\sharp}$ uniformly~$\nu \geq 1$.
  Clearly, both these ingredients hold in the space~$\sL^p_{\nu,
    \sharp}$, $1 \leq p \leq \infty$.  Thus, the Haar
  system~\eqref{Inductive}, is a basis in~$\sL^p_{\nu, \sharp}$.  More
  precisely, the constant~$c^{(p)}_{\nu, \sharp}$ which guarantee the
  inequality
    $$
    \|\sum_{0\le j \leq m} \alpha_j h^{(\nu)}_j \|_{p, \sharp} \le
    c^{(p)}_{\nu, \sharp}\, \|\sum_{0\le j < 4^\nu} \alpha_j h
    ^{(\nu)}_j \|_{p, \sharp}, $$
    for every~$0\le m < 4^\nu$ and every
    complex scalars~$\alpha_j$, admits the estimate $$
    c^{(p)}_{1,
      \sharp} \le \sum_{0\le j\le 3} \|r_j\| \, \|r_j\|_{p,
      \sharp}
  $$
  and $$
  c^{(p)}_{\nu +1, \sharp} \le \max\{
  c^{(p)}_{\nu, \sharp}\, \|r_0\|^2; \|r_0\|^2 +
  2\, \sum_{q=1}^3 \|r_q\|^2 \}. $$
  The similar estimates hold true for the right spaces~$\sL^p_{\nu,
    \flat}$ and the right Haar system.
\end{rem}

\begin{rem}
  The construction of the Haar system may be generalized as follows: in
  the inductive definition~\eqref{Inductive} for each inductive step
  from~$\h_\nu$ to~$\h_{\nu+1}$ we can use its own set~$\{r^{(\nu)}_0,
  r^{(\nu)}_1, r^{(\nu)}_2, r^{(\nu)}_3\}$ which possesses the
  property~\eqref{basis} (or~\eqref{leftBasis}, if we build a
  right Haar system).  Theorem~\ref{MTh} remains valid in this
  case with obvious changes to the estimates of the constants~$c_{\nu,
  \sharp (\text{resp.}\  \flat)}$.
\end{rem}

We shall refer to the system~$\h_\nu$ constructed above
as~$\h^{(\nu)}_\alpha(r_0, r_1, r_2, r_3)$ in the sequel to stress the
fact that the Haar system~$\h^{(\nu)}_\alpha$ depends on~$0 < \alpha
\leq \frac 12$ and~$\{r_j\}_{0 \leq j \leq 3}$.

\begin{corl}
  \label{unifMth}
  The system~$\h^{(\nu)}_\alpha (r_0, r_1, r_2, r_3)$,
  where~$\{r_j\}_{0 \leq j \leq 3}$ satisfies~\eqref{basis}
  (resp.~\eqref{leftBasis}) and~$\|r_0\|\leq 1$, is a basis
  in~$\aN_{\nu, \sharp}$ (resp.~$\aN_{\nu, \flat}$) with the basis
  constant uniformly bounded with respect to~$\nu \in \N$.
\end{corl}

\subsubsection{The hyperfinite factors~$\hbox{II}_1$
  and~$\hbox{III}_\lambda$, $0< \lambda < 1$.}

The collection of the algebras~$\{(\aN_\nu, \rho_\nu)\}_{\nu\in\N}$
together with the embedding~\eqref{TowerEmbed} forms {\it a directed
  system of $C^*$-algebras}, \cite[Section~11.4]{KadRingII}.  {\it The
  inductive limit\/} of this system possesses a state~$\rho_\alpha$,
induced by~$\rho_\nu$, $\nu \geq 1$.  We denote the GNS representation
of this inductive limit with respect to the state~$\rho_\alpha$
as~$\HFF_\alpha$.  $\HFF_\alpha$ is a factor of
type~$\hbox{III}_\lambda$ if~$0< \alpha < \frac 12$, with~$\lambda =
\frac \alpha{1 - \alpha}$ and a factor of type~$\hbox{II}_1$ if~$\alpha
= \frac 12$.  The properties of the factor~$\HFF_\alpha$ are collected
in the following lemma.  We also refer the reader
to~\cite[Section~12.3]{KadRingII}, where the representation of the
factor~$\HFF_\alpha$ as a discrete crossed product is given.

\begin{lem}
  \label{matricial}
  $\HFF_\alpha$ is a matricial factor which possesses a distinguished
  faithful normal state~$\rho_\alpha$.  With~$\aN_\nu$, $\rho_\nu$,
  $\nu\in\N$ defined in the previous section, there are ultra-weakly
  continuous $*$-isomorphic embeddings $\pi_\nu:\aN_\nu \rightarrow
  \HFF_\alpha$, $\nu\in\N$ such that
  \begin{enumerate}
  \item the embedding~$i_\nu: \aN_\nu \rightarrow \aN_{\nu+1}$, given
    in~\eqref{TowerEmbed}, carries into~$\pi_\nu(\aN_\nu) \subseteq
    \pi_{\nu+1}(\aN_{\nu+1})$;
  \item the state~$\rho_\nu$ is induced by~$\rho_\alpha$ and~$\pi_\nu$,
    i.e.\ $\rho_\nu(x) = \rho_\alpha(\pi_\nu(x))$, $x\in\aN_\nu$,
    moreover, the automorphism group~$\sigma^{\rho_\alpha}$ leaves
    every subalgebra~$\pi_\nu(\aN_\nu)$ globally invariant;
  \item the set~$\cup_{\nu \geq 1} \pi_\nu(\aN_\nu)$ is dense
    in~$\HFF_\alpha$ with respect to the weak operator topology.
  \end{enumerate}
\end{lem}

From now on we shall identify the algebras~$\aN_\nu$
with~$\pi_\nu(\aN_\nu)$, $\nu\in\N$.  Since the
group~$\sigma^{\rho_\alpha}$ leave the subalgebras~$\aN_\nu$, $\nu\in\N$
globally invariant, it follows from the preliminaries that, for
every~$1\le p \le \infty$, we have
\begin{equation}
  \label{PredualTower}
  \sL^p_{1, \sharp} \subseteq \ldots \subseteq
  \sL^p_{\nu, \sharp} \subseteq \sL^p_{\nu+1, \sharp} \subseteq 
  \ldots \subseteq L^p_{\sharp}(\HFF_\alpha),
\end{equation}
and all embeddings here are isometric.  Thus, we may refer to the norms
in all these spaces as~$\|\cdot\|_{p,\sharp}$, omitting the index~$\nu\in
\N$.  Moreover, these embeddings are $1$-complemented, i.e.\ there is norm
one projection
\begin{equation}
\label{CERlambda}
\CE_\nu: L^p_{\sharp}(\HFF_\alpha) \rightarrow
\sL^p_{\nu, \sharp},\ \  \nu \in\N. 
\end{equation}
We also note that~$\cup_{\nu \in\N} \sL^p_{\nu, \sharp}$ is norm dense
in~$L^p_{\sharp}(\HFF_\alpha)$, $1\le p < \infty$, since
\begin{equation}
  \label{Convergence}
  \lim_{\nu \rightarrow \infty} \|\CE_\nu(x) - x \|_{p, \sharp} = 0,\ \ 
  x\in L^p_{\sharp}(\HFF_\lambda),\ \ 1\le p < \infty.
\end{equation}
The last statement (and its right counterpart) is established
in~\cite[Theorem~8]{Go1985}.

Since the space~$\HFF_\alpha$ is not separable, the
convergence~\eqref{Convergence} can not be extended to the norm~$\|\cdot\|$.
Nonetheless, we have the ultra-weak convergence in this case, namely

\begin{lem}
  \label{WeakConv}
  For every~$x \in \HFF_\alpha$ and~$\phi\in \HFF_{\alpha, *}$, we have
  \begin{equation*}
    \lim_{\nu \rightarrow \infty} \phi(\CE_\nu(x) - x) = 0.
  \end{equation*}
\end{lem}

\begin{proof}
  The proof is straitforward.  At first, from~\eqref{Convergence}, we
  have
  $$ |\rho_\alpha(y\, (\CE_\nu(x) - x))| \le \, \|y\| \, \|\CE_\nu(x) -
  x\|_\sharp \rightarrow 0,\ \text{as}\ \nu \rightarrow
  \infty,\ \ x, y\in \HFF_\alpha.  $$
  The latter means that we proved the lemma for the special case~$\phi
  = \rho y$, $y\in\HFF_\alpha$.  Since the linear subspace~$\{\rho
  y\}_{y \in \HFF_\alpha}$ is norm dense in~$\HFF_{\alpha, *}$ and the
  projections~$\CE_\nu :\HFF_\alpha \rightarrow \aN_\nu$ are uniformly
  bounded, $\nu \in \N$, the general case now follows
  from~\cite[Lemma~1.2]{SZ-Lec-on-vNA}.
\end{proof}

Alternatively, we may look at the tower~\eqref{PredualTower} from {\it
  inductive limit\/} point of view as follows
(see~\cite[p.~135]{Pisier-Factorization}) for definition of inductive
limits of Banach spaces)

\begin{thm}
  \label{InLimI}
  If~$1\le p < \infty$, the collection of Banach spaces~$\{\sL^p_{\nu,
    \sharp}\}_{\nu\in\N}$ together with the
  embedding~\eqref{TowerEmbed} forms {\it a directed system of Banach
    spaces.}  The inductive limit of this system is isomorphic
  to~$\sL^p_{\sharp}(\HFF_{\alpha})$.
\end{thm}

Let us further note that the identities~\eqref{leftRightPreNorms}
and~\eqref{LpLeftRightIdent} mean that the Banach spaces~$\sL^p_{\nu,
  \sharp}$ and~$\cC_p^{(\nu)}$ are isometric with the isometry given by
\begin{equation}
  \label{NTraceIso}
  x \in \sL^p_{\nu, \sharp} \rightarrow x
  \Alpha^{1/p}_\nu \in \cC_p^{(\nu)}.
\end{equation}
Applying this isometry to the directed system~$\{\sL^p_{\nu,
  \sharp}\}_{\nu \in \N}$ together with the embedding~$i_\nu$ we obtain
\begin{corl}
  \label{InLimII}
  The collection of matrix spaces~$\{\cC_p^{(\nu)}\}_{\nu\in\N}$, $1\le
  p < \infty$, equipped with the $p$-th Schatten-von Neumann norm,
  together with the embedding
  $$
  x\in \cC_p^{(\nu)} \rightarrow x \otimes \Alpha^{1/p}_1 \in
  \cC^{(\nu+1)}_p $$
  is a directed system of Banach spaces with the
  inductive limit isomorphic to Haagerup's space~$L^p(\HFF_{\alpha})$.
\end{corl}

\subsubsection{Haar system (cont.)}

Consider the left (resp.\ right) Haar system~$\h^{(\nu)}_\alpha =
\h^{(\nu)}_\alpha (r_0, r_1, r_2, r_3)$, where the system~$\{r_j\}_{0
  \leq j\leq 3}$ satisfies~\eqref{basis} (resp.~\eqref{leftBasis}) such
that~$r_0 = \id_1$.  It then follows from~\eqref{TowerEmbed}
and~\eqref{Inductive} that~$h^{(\nu+1)}_j = i_\nu(h^{(\nu)}_j)$, $0\le j
< 4^\nu$, $\nu\in\N$.  Thus, we can construct a unified left (resp.\ 
right) Haar system~$\h_\alpha = \h_\alpha(r_0, r_1, r_2, r_3) =
\{h_j\}_{j\ge 0}$ in~$\HFF_\alpha$ as~$h_j = \pi_\nu(h_j^{(\nu)})$,
provided~$0\le j < 4^\nu$.  As a corollary of Theorem~\ref{MTh}, we now
have

\begin{thm}
  \label{IIIlambdaMth}
  The left (resp.\ right) Haar system~$\h_\alpha$ forms a basis in the
  space~$\HFF_{\alpha, \sharp (\text{resp.}\  \flat)}$.
\end{thm}

\begin{proof}
  To prove that~$\h_\alpha$ is a basis, we need to prove
  (i)~$\h_\alpha$ is a basic sequence and (ii)~the linear span
  of~$\h_\alpha$ is dense in~$\HFF_{\alpha, \sharp}$ (resp.\
  $\HFF_{\alpha, \flat}$).  The first part is contained in
  Theorem~\ref{MTh} and the second one is guaranteed
  by~\eqref{Convergence}.
\end{proof}

As an example of the system~$\h^{(1)}_\alpha$, we now take the system
\begin{equation}
  \label{typicalSystem}
  r_0 = \id_1,\ 
  r_1 = \left[\begin{matrix} \frac 1 {\sqrt\lambda}  & 0 \\ 0 & -
      \sqrt\lambda \\
    \end{matrix}\right],\
  r_2 = \left[\begin{matrix} 0 & 1 \\ 1 & 0 \\ \end{matrix}\right],\ 
  r_3 = \left[\begin{matrix} 0 & - \sqrt \lambda \\ \frac 1 
      {\sqrt \lambda} & 0 
    \end{matrix}\right],
\end{equation}
where~$\lambda = \frac\alpha{1- \alpha}$.  The system~$\h^{(1)}_\alpha$
satisfies~\eqref{basis} (resp.~\eqref{leftBasis}, if~$r_3$ replaced
with~$r^*_3$).  Thus, from Theorem~\ref{IIIlambdaMth}, we have

\begin{corl}
  \label{typicalTh}
  The left (resp.\ right) Haar system~$\h_\alpha = \h_\alpha(r_0, r_1,
  r_2, r_3)$, where~$\{r_j\}_{0 \leq j \leq 3}$ are given
  in~\eqref{typicalSystem}, is a basis in the predual of the hyperfinite
  factor~$\hbox{III}_\lambda$.
\end{corl}

For the special case~$\alpha=\frac 12$, the system~\eqref{typicalSystem}
turns into
\begin{equation}
  \label{typicalSystemI}
  \hat r_0 = \id_1,\ 
  \hat r_1 = \left[\begin{matrix} 1  & 0 \\ 0 & - 1 \end{matrix}\right],\
  \hat r_2 = \left[\begin{matrix} 0 & 1 \\ 1 & 0 \\ \end{matrix}\right],\ 
  \hat r_3 = \left[\begin{matrix} 0 & - 1 \\ 1 & 0 \end{matrix}\right],
\end{equation}
and we also have

\begin{corl}
  \label{typThI}
  The Haar system~$\h_{\frac 12} = \h_{\frac 12}(\hat r_0, \hat r_1,
  \hat r_2, \hat r_3)$ is a basis in the predual of the hyperfinite
  factor~$\hbox{II}_1$~$\HFF_{\frac 12}$.
\end{corl}

Let us next consider the diagonal subalgebras~$\aA_\nu \subseteq
\aN_\nu$, $\nu\in\N$.  The weak-operator closure of~$\cup_{\nu\in\N}
\pi_\nu(\aA_\nu)$ forms an Abelian subalgebra~$\aA_\alpha$
in~$\HFF_\alpha$, which is isomorphic to~$L_\infty([0, 1), m_\alpha)$,
cf.~\cite[Section~12.3]{KadRingII}, the algebra of all essentially
bounded $m_\alpha$-measurable functions on~$[0,1)$, where the
measure~$m_\alpha$ is given by
\begin{equation}
  \label{DistortedMeasure}
  m_\alpha([\frac k{2^\nu}, \frac {k+1}{2^\nu}]) = \prod_{s = 0}^{\nu-1}
  [(1-\epsilon_s) \alpha + \epsilon_s (1-\alpha)],
\end{equation}
where~$0\le k < 2^\nu$ and~$\epsilon_s$ are binary digits of~$k$,
i.e.~$\epsilon_s = 0, 1$ such that $$ k = \epsilon_0 2^0 + \epsilon_1
2^1 + \ldots + \epsilon_{\nu-1} 2^{\nu-1}. $$ Since~$\aA_\alpha$ is
commutative, we have that~$\aA_{\alpha, \sharp} = \aA_{\alpha, \flat}
= \aA_{\alpha, *}$.  The modular automorphism
group~$\sigma^{\rho_\alpha}$ (resp.\ $\sigma^{\rho_\nu}$) leaves the
subalgebra~$\aA_\alpha$ (resp.\ $\aA_\nu$) globally invariant.  Thus,
the embedding~$\aA_{\alpha, *} \subseteq \HFF_{\alpha, \sharp
  (\text{resp.}\ \flat)}$ (resp.\ $\aA_{\nu, *} \subseteq \aN_{\nu,
  \sharp (\text{resp.}\ \flat)}$) is isometric and complemented.
Let~$\CE$ be the norm one projection~$\CE:\HFF_\alpha \rightarrow
\aA_\alpha$.  We denote by the same letter~$\CE$ the norm one
projection~$\CE: \HFF_{\alpha, \sharp (\text{resp.}\ \flat)}
\rightarrow \aA_{\alpha, *}$.  The projection~$\CE:\aN_\nu \rightarrow
\aA_\nu$ vanishes on all non-diagonal matrix entries.  Hence, we
obtain that, if~$\h_\alpha = \h_\alpha(r_0, r_1, r_2, r_3) =
\{h_j\}_{j\ge 0}$ is the left (or right) Haar system, with respect
to~\eqref{typicalSystem}, then $$ \CE(h_j) = h_j,\ \ \
\vtop{\hsize=0.6\hsize \noindent if $4^\nu \le j < 2\cdot 4^\nu$
  and~$h_j = i_\nu(e_{k}^{(\nu)})\cdot (\id_\nu \otimes r_1)$ with~$k
  = j - 4^\nu$ and~$e^{(\nu)}_k$ being a diagonal matrix unit,
  see~\eqref{Inductive};} $$
\begin{equation}
  \label{commutativeHaar}
  \CE(h_j) = 0,\ \ \ \vtop{\hsize=0.6\hsize \noindent otherwise.}
\end{equation}
Clearly, this implies that the non-zero subsystem~$\chi_\alpha =
\{\chi_j\}_{j\ge0}$ of~$\CE(\h_\alpha)$ form a basis of~$\aA_{\alpha,
  *}$.  From~\eqref{Inductive} and~\eqref{commutativeHaar}, we obtain
that, $\chi_0 = r_0$, $\chi_1 = r_1$ and, if~$2^\nu \le j <
2^{\nu+1}$,
\begin{equation}
  \label{CommutativeInductive}
  \chi_j = \varepsilon^{(\nu)}_k \cdot (\id_\nu \otimes r_1),\ \ 
  j = 2^\nu + k,\ \ 0\le k < 2^\nu,
\end{equation}
where~$\varepsilon^{(\nu)}_k$ is the $k$-th diagonal matrix unit
in~$\aN_\nu$.  When~$\alpha = \frac 12$, the system~$\chi_{\frac 12}$
is the classical Haar system, cf.~\cite[Section~2.c]{LT-II}.  Thus, we
have

\begin{corl}
  \label{CommutativeMTh} The system~$\chi_\alpha$, given
  in~\eqref{CommutativeInductive} is a basis of~$L_1([0,1), m_\alpha)$.
  In particular, for~$\alpha = \frac 12$ this system coincides with the
  classical Haar system on~$L_1(0, 1)$, \cite[Section~2.c]{LT-II}.
\end{corl}

\begin{rem}
  \label{LpAgain}
  Every result in this sections extends to the $L^p$-spaces associated
  with the factors~$\hbox{III}_\lambda$ and~$\hbox{II}_1$.  If~$p =
  \infty$, then the results still hold true with norm convergence
  replaced by ultra-weak convergence, cf.\ Lemma~\ref{WeakConv}.
\end{rem}

\begin{rem}
  \label{Rademachers} In analogy with the classical Haar system, we
  shall call the system~\eqref{typicalSystem}, and its
  derivatives~$r_{j, \nu} = \id_\nu \otimes r_j$, $0 \leq j\leq 3$,
  $\nu \in \N$ {\it the (non-commutative) Rademacher system.}  Due to
  unconditionality of martingale differences in the
  spaces~$L^p(\HFF_{\frac 12})$, $1 < p < \infty$, \cite{SF1994,
    SF1995, PiXu1997}, the Rademacher system is an unconditional basic
  sequence in~$L^p(\HFF_{\frac 12})$, $1< p< \infty$.
\end{rem}

\subsubsection{Factors of type~$\hbox{III}_1$ and~$\hbox{II}_\infty$.}

Here we shall consider the construction of bases in the preduals of
the factors of type~$\hbox{III}_1$ and~$\hbox{II}_\infty$.  Since
these two factors may be reduced to the factors of
type~$\hbox{III}_\lambda$, $\hbox{II}_1$ and~$\hbox{I}_\infty$ by
means of tensor products, we shall first consider the extension of the
Haar system construction over preduals of tensor products.  To this
end, it is useful to recall the notion of {\it Schauder
  decomposition}, \cite{LT-I}.

Let~$\D = \{D_j\}_{j \geq 1}$ be a system of projections in a Banach
space~$X$ such that~$D_j D_k = 0$, $j \neq k$.  The system~$\D$ is a
Schauder decomposition of the Banach space~$X$ if and only if the
series~$\sum_{j =1}^\infty D_j x$ converges to~$x$ in the norm of~$X$,
for every~$x \in X$.  As for bases, we have the equivalent criteria
for the system~$\D$ to be a Schauder decomposition of~$X$,
\cite{LT-I}.  Namely, a system~$\D = \{D_j\}_{j \geq 1}$, $D_j D_k =
0$, $j \neq k$ is a Schauder decomposition of~$X$ if and only if
(i)~$[D_j(X)]_{j \geq 1} = X$; (ii)~there is a constant~$c$ such that
$$ \|\sum_{1 \leq j \leq m} D_j x \|_X \leq c\, \|\sum_{ 1\leq j \leq
  n} D_j x \|_X,\ \ x \in X, \ 1\leq m \leq n. $$

We fix two von Neumann algebras~$\aN$ and~$\aM$ equipped with faithful
normal states~$\rho$ and~$\phi$, respectively.  The tensor product
algebra~$\aN \bar \otimes \aM$ with respect to the product state~$\rho
\otimes \phi$ is the weak operator completion of the GNS
representation of the tensor product $C^*$-algebra~$\aN \otimes \aM$
with respect to the state~$\rho \otimes \phi$,
cf.~\cite[Chapter~11]{KadRingII}.

We shall consider the algebra~$\aN$ as von Neumann subalgebra of~$\aN
\otimes \aM$ under the embedding~$x \rightarrow x \otimes \id$, $x\in\aN$.
It clearly follows from the modular condition that $$ \sigma_t^{\rho
  \otimes \phi}(x \otimes y) = \sigma_t^{\rho}(x) \otimes
\sigma_t^{\phi}(y),\ \ x\in \aN, y \in \aM, t \in \Rl. $$ Thus, the
modular group~$\sigma^{\rho \otimes \phi}$ leaves subalgebra~$\aN$
globally invariant, and therefore, the results in the preliminaries
are applicable.  In particular, the left predual~$\aN_{\sharp}$ is
isometrically embeds into~$(\aN \bar \otimes \aM)_{\sharp}$, and the
space~$\aN_{\sharp}$ is $1$-complemented in~$(\aN \bar \otimes
\aM)_\sharp$.  Let us denote the corresponding projection
as~$\CE_\aN$.  As in~\eqref{CEprop}, we obtain the explicit formula
for the projection~$\CE_\aN$ on elementary tensors
\begin{equation}
  \label{CEpropI}
  \CE_\aN(x \otimes y) = (x \otimes \id)\,
  \phi(y),\ \  x\in\aN,\ y \in \aM.
\end{equation}

Let us fix an orthonormal basis~$\y = \{y_j\}_{j \geq 1} \subseteq
\aM$ in the predual~$\aM_{\sharp}$, $i=1,2$.  We assume that
\begin{equation}
  \label{DualSys}
  \phi (y^*_j y_k) = \delta_{jk}. 
\end{equation}
Having basis~$\y$ and the expectation~$\CE_\aN$ at our disposal, we can
construct the associated system of projections~$\D = \{D_j\}_{j \geq
  1}$ of the predual~$(\aN \bar \otimes \aM)_\sharp$ by
\begin{equation}
  \label{AccosDecomp}
  D_j z = (\id \otimes y_j)\, \CE_\aN((\id \otimes
  y_j)^*\, z),\ \ z\in \aN \bar \otimes \aM.
\end{equation}
The left multiplication by an element of the tensor product~$\aN \bar
\otimes \aM$ is a bounded operator on~$(\aN \bar \otimes \aM)_\sharp$.
Therefore, the operators~$D_j$ are indeed bounded linear operators
on~$(\aN \bar \otimes \aM)_\sharp$.  The fact that~$D_j$ are projections
and that~$D_j D_k = 0$ if~$j \neq k$ follows from~\eqref{CEpropI}
and~\eqref{DualSys}.  Furthermore, the identities~\eqref{CEpropI}
and~\eqref{DualSys} give the explicit formula for the projection~$D_j$
on the algebraic tensor product~$\aN \otimes \aM$.  Indeed, if~$a_k \in
\aN$, $b_k \in \aM$, $1 \leq k \leq n$ and $$
b_k = \sum_{s \geq 1}
\alpha_{ks}\, y_s,\ \ \alpha_{ks} \in \Cx $$
is the expansion of the
element~$b_k$ with respect to the system~$\y$ in~$\aM_{\sharp}$, then
\begin{equation}
  \label{AlgTensFor}
  D_j (\sum_{1 \leq k \leq n} a_k \otimes b_k) = \sum_{1 \leq k \leq n}
  \alpha_{kj} \, a_k \otimes y_j. 
\end{equation}
Since the algebraic tensor product~$\aN \otimes \aM$ is norm dense
in~$(\aN \bar \otimes \aM)_\sharp$ and the norm in the space~$(\aN
\bar \otimes \aM)_\sharp$ is a cross-norm, we get~$[D_j(\aN \bar
\otimes \aM)]_{j \geq 1} = (\aN \bar \otimes \aM)_\sharp$.  In
general, it is not the case that the system~$\D$ constructed above is
a Schauder decomposition of the Banach space~$(\aN \bar \otimes
\aM)_\sharp$.

\begin{thm}
  \label{TesorBasis} If~$\x = \{x_j\}_{j \geq 1}$ and~$\y = \{y_k\}_{k
    \geq 1}$ are two bases in~$\aN_\sharp$ and~$\aM_\sharp$,
  respectively, such that the associated systems of projections~$\E =
  \{E_j\}_{j \geq 1}$ and~$\D = \{D_k\}_{k \geq 1}$, defined
  in~\eqref{AccosDecomp}, are Schauder decompositions of the Banach
  space~$(\aN \bar \otimes \aM)_\sharp$, then the product basis~$\z =
  \x \otimes \y := \{x_j \otimes y_k\}_{j, k \geq 1}$, taken in the
  shell enumeration, is a basis in the predual~$(\aN \bar \otimes
  \aM)_\sharp$.
\end{thm}

\begin{proof}
  The proof is rather standard.  The shell enumeration assigns to a
  pair~$(j, k)$, $j,  k \geq 1$ the number~$s(j,k)$ defined by
  $$  s(j, k) = \begin{cases} (k-1)^2 + j, & \text{if~$j \leq k$;} \\
    j^2 - k + 1, & \text{if~$j > k$.} \\ \end{cases} $$ Let~$\z =
  \{z_s\}_{s \geq 1}$ be the system~$\z$ in the shell enumeration.  It
  is clear that the linear span of the system~$\z$ is dense in~$(\aN
  \bar \otimes \aM)_\sharp$.  Thus, to prove the theorem, we have to
  establish, that there is a constant~$c$ such that $$ \|\sum_{1 \leq
    s \leq m} \alpha_s z_s\|_\sharp \le c\, \|\sum_{1 \leq s \leq n}
  \alpha_s z_s \|_\sharp,\ \ \alpha_s\in \Cx, 1\leq m \leq n. $$
  Without loss of generality we may assume that~$n = n_1^2$, for
  some~$n_1 \geq 1$.  There is an integer~$m_1 \geq 1$ such that
  either of the relations is true (i)~$m_1^2+1 \leq m \leq m_1^2 +
  m_1$ or (ii)~$m_1^2+m_1+1 \leq m \leq (m_1+1)^2$.  Let us consider
  the first option, for the second one the argument is similar.
  We have
  \begin{equation}
    \label{temp}
    \begin{aligned}
      \|\sum_{1\leq s \leq m} \alpha_s z_s \|_\sharp \leq & \, \|\sum_{1
        \leq s \leq m_1^2 } \alpha_s z_s \|_\sharp + \|\sum_{m_1^2 + 1
        \leq s \leq m} \alpha_s z_s \|_\sharp \\ = & \, \|\sum_{1\leq
        j,k \leq m_1} \alpha_{s(j,k)}\,
      x_j \otimes y_k \|_\sharp \\ + &\, \| \sum_{1 \leq j
        \leq m - m_1^2 - 1} \alpha_{s(j, m_1)}\, x_j \otimes
      y_{m_1} \|_\sharp
    \end{aligned}
  \end{equation}
  Letting $$ z = \sum_{1 \leq s \leq n} \alpha_s z_s = \sum_{1 \leq
    j,k \leq n_1} \alpha_{s(j,k)} x_j \otimes y_k, $$ it then follows
  from~\eqref{AlgTensFor} that latter two terms on the right hand side
  in~\eqref{temp} are given by $$ \sum_{1\leq j,k \leq m_1}
  \alpha_{s(j,k)}\, x_j \otimes y_k = P_{m_1} Q_{m_1}(z) $$ and $$
  \sum_{1 \leq j \leq m - m_1^2 - 1} \alpha_{s(j, m_1)}\, x_j \otimes
  y_{m_1} = P_{m - m_1^2 -1} D_{m_1} (z), $$ where~$P_j$ and~$Q_k$ are
  partial sum projections with respect to the decompositions~$\E$
  and~$\D$, respectively, i.e. $P_j = \sum_{1 \leq s \leq j} E_s$
  and~$Q_k = \sum_{1 \leq s \leq k} D_s$.  Thus, we continue
  \begin{equation*}
    \|\sum_{1\leq s \leq m} \alpha_s z_s \|_\sharp \leq \|P_{m_1}
    Q_{m_1}(z)\|_\sharp + \|P_{m - m_1^2 - 1} D_{m_1}(z)\|_\sharp 
     \leq c \, \|z\|_\sharp.
  \end{equation*}
  The latter inequality is due to the fact that the partial sum
  projections are uniformly bounded.  The claim of the theorem follows.
\end{proof}

We consider two specific examples.  Firstly, we take~$\aM =
\Bd(\ell_n^2)$.  In this setting the algebra~$\aN \otimes \aM$ may be
considered as the space of all bounded $n \times n$-block matrices
with entries in~$\aN$.  The latter observation means, in particular,
that the system of projections~$\E = \{E_j\}_{j \geq 1}$,
\eqref{AccosDecomp}, associated with the matrix unit basis~$\e$ is a
Schauder decomposition of the predual~$(\aN \bar \otimes \aM)_\sharp$,
see the proof Lemma~\ref{MatrixUnits}.  Let us also note that the
Schauder constant of the system~$\E$ is uniformly bounded with respect
to dimension of~$\ell^2_n$.

As another example, we take~$\aM = \HFF_\alpha$, $0 < \alpha \leq
\frac 12$.  We fix the Haar system~$\h_\alpha = \h_\alpha(r_0, r_1,
r_2, r_3) = \{h_j\}_{j \geq 1}$, where the Rademachers~$\{r_j\}_{0
  \leq j \leq 3}$ are given in~\eqref{typicalSystem}.  We denote the
associated system of projections by~$\H_\alpha = \H_\alpha(r_0, r_1,
r_2, r_3) = \{H_j\}_{j \geq 1}$.  To establish, that the
system~$\H_\alpha$ is a Schauder decomposition in~$(\aN \bar \otimes
\aM)_\sharp$ we only need to verify, that there is a constant~$c$ such
that $$ \|\sum_{1 \leq j \leq m} H_j z \|_\sharp \le c \, \| \sum_{1
  \leq j \leq n} H_j z \|_\sharp,\ \ z\in \aN \otimes \aM,\ \ 1\leq m
\leq n. $$ The proof of the latter inequality is based on the
following theorem.

\begin{thm}
  \label{MThDecom}
  Let~$(\aN_\nu, \rho_\nu)$ be the algebras from the preceding sections
  and let\/~$\aN \bar \otimes \aN_\nu$ be the tensor product von Neumann
  algebras equipped with the product states~$\rho \otimes \rho_\nu$,
  $\nu \in \N$.  Let\/~$\h^{(\nu)}_\alpha = \h^{(\nu)}_\alpha (r_0, r_1,
  r_2, r_3) = \{h_j^{(\nu)}\}_{j \geq 1}$, where~$\{r_j\}_{0 \leq j \leq
    3}$ satisfies~\eqref{basis}, be the left Haar system in~$\aN_\nu$.
  If\/~$\H^{(\nu)}_\alpha = \H^{(\nu)}_\alpha(r_0, r_1, r_2, r_3)$ is the
  associated decomposition defined by~\eqref{AccosDecomp}, then the
  minimal constant~$c$, which guarantees the inequality $$
  \|\sum_{0
    \leq j < m} H_j^{(\nu)} z \|_\sharp \le c_{\nu, \sharp}\, \|\sum_{0
    \leq j < 4^{\nu}} H^{(\nu)}_j z\|_\sharp,\ \ z \in \aN \bar \otimes
  \aN_\nu,\ 0 \leq m < 4^\nu, $$
  admits the same inductive estimate as
  that in Theorem~\ref{MTh}.
\end{thm}

\noindent The proof of Theorem~\ref{MThDecom} is essentially a
repetition of that of Theorem~\ref{MTh}, we leave details to the
reader.  Thus, we obtain

\begin{thm}
  \label{SchDecRes}
  The system~$\H_\alpha = \H_\alpha(r_0, r_1, r_2, r_3) = \{H_j\}_{j
    \geq 1}$, associated with the Haar basis~$\h_\alpha(r_0, r_1, r_2,
  r_3)$, where~$\{r_j\}_{0 \leq j \leq 3}$ given
  by~\eqref{typicalSystem}, is a Schauder decomposition of the Banach
  space~$(\aN \bar \otimes \HFF_\alpha)_\sharp$.
\end{thm}

Now we may apply the results above to the hyperfinite
factors~$\hbox{III}_1$ and~$\hbox{II}_\infty$.  It is known
that~$\hbox{III}_1 = \hbox{III}_{\lambda_1} \bar \otimes
\hbox{III}_{\lambda_2}$, where~$\frac {\log \lambda_1 }{\log
  \lambda_2}\notin \Q$ and~$\hbox{II}_\infty = \hbox{II}_1 \bar
\otimes \hbox{I}_\infty$, cf.~\cite{Co1976, Ha1987}.  See the
definition of~$\h_\alpha$ and~$\h_{\frac 12}$ in
Corollaries~\ref{typicalTh} and~\ref{typThI}.
obtain

\begin{corl}
  \label{III1base}
  The system~$\z = \h_{\alpha_1} \otimes \h_{\alpha_2}$, $\alpha_i =
  \frac {\lambda_i}{\lambda_i + 1}$, $i=1,2$, $\frac{\log
    \lambda_1}{\log \lambda_2} \notin \Q$, is a basis in the predual
  of the hyperfinite factor of type~$\hbox{III}_1$.
\end{corl}

\begin{corl}
  \label{IIinfbase}
  The system~$\x = \h_{\frac 12} \otimes \e$, is a basis in the
  predual of the hyperfinite factor of type~$\hbox{II}_\infty$.
\end{corl}

\begin{rem}
  Similarly to preceding sections, all results above hold true in the
  setting of left and right $L^p$-spaces associated with the factors
  of type~$\hbox{III}_1$ and~$\hbox{II}_\infty$.  Moreover, the
  system~$h_{\frac 12}$ (resp.~$\h_{\frac 12} \otimes \e$) from
  Corollary~\ref{typThI} (resp.~\ref{IIinfbase}) forms a basis in any
  symmetric operator space~$E(\HFF_{\frac 12})$ (resp.~$E(\aM)$),
  where~$\HFF_{\frac 12}$ (resp.~$\aM$) is a hyperfinite
  factor~$\hbox{II}_1$ (resp.~$\hbox{II}_\infty$) and~$E$ is separable
  rearrangement invariant function space (see definitions and further
  references in~\cite{DFPS2001}).
\end{rem}

{\baselineskip=0.8\baselineskip \small \parskip=0.8\parskip
  \let\Large=\normalsize \bibliography{references}}

\medskip

\hrule width 5cm

\medskip

\parindent=0pt

\small

D.~Potapov, pota0002@infoeng.flinders.edu.au, \\
Dr. F.~Sukochev, sukochev@infoeng.flinders.edu.au,

\medskip

School of Informatics and Engineering, \\
Faculty of Science and Engineering, \\
Flinders Univ. of SA, Bedford Park, 5042, \\
Adelaide, SA, Australia.

\end{document}